\documentclass{amsart}

\usepackage[lmargin=3.5cm,rmargin=3.5cm,top=4cm,bottom=4cm]{geometry}
\newtheorem{theorem}{Theorem}[section]
\newtheorem{proposition}[theorem]{Proposition}

\newtheorem{corollary}[theorem]{Corollary}
\newtheorem{definition}[theorem]{Definition}
\newtheorem{remark}[theorem]{Remark}

\numberwithin{equation}{section}
\usepackage{caption} 
\usepackage{subcaption} 
\usepackage{graphicx}
\usepackage{pgfplots}
\usepackage[all]{nowidow}
\usepackage[utf8]{inputenc}
\usepackage{tikz}
\usepackage[inline]{enumitem}
\usepackage{blindtext}
\usepackage{amssymb}		
\usepackage{ulem}
\usepackage{amsfonts}

\usepackage[hidelinks,colorlinks=true,linkcolor=blue,citecolor=red]{hyperref}
\numberwithin{equation}{section}

\begin{document}

\title{A Variational Method to Calculate Probabilities}


\author[Hugo Reyna Castañeda]{Hugo Guadalupe Reyna-Castañeda}
\address{Departamento de Matemáticas, Facultad de Ciencias, Universidad Nacional Autónoma de México, Mexico City, Mexico}
\curraddr{}
\email{hugoreyna46@ciencias.unam.mx}
\thanks{}

\author{Mar\'ia de los \'Angeles Sandoval-Romero}
\address{Departamento de Matemáticas, Facultad de Ciencias, Universidad Nacional Autónoma de México, Mexico City, Mexico}
\curraddr{}
\email{selegna@ciencias.unam.mx}
\thanks{}


\date{}

\dedicatory{}

\begin{abstract}
In this paper, we prove the existence and uniqueness of the conditional expectation of an event $A$ given a $\sigma$-algebra $\mathcal{G}$ as a linear problem in Lebesgue spaces $L^{p}$ associated with a probability space through the Riesz Representation Theorems. For the $L^{2}$ case, we state Dirichlet's principle. Then, we extend this principle to specific values of $p$, defining the existence of the conditional expectation as a variational problem. We conclude with a proof of the law of total probability using these tools.
\end{abstract}

\maketitle

\tableofcontents

\section{Introduction}

The law of total probability (see \cite[Theorem 3.4]{Jacod(2004)}) is a fundamental result of classical probability theory with multiple applications in various disciplines such as physics, medicine, biology, etc. This theorem allows for the calculation of an event probability through the computation of conditional probabilities.

\begin{theorem}[Law of total probability]\label{TotalProbability}
Let $(\Omega,\mathcal{F},P)$ be a probability space and let $\mathcal{B}=\{B_1,\ldots,B_N\}$ be a collection on $\mathcal{F}$ such that $\Omega=\bigcup_{j=1}^{N} B_j$, $B_{i} \cap B_{j} = \varnothing$ for all $i \neq j$ and $P(B_j) >0$ for all $j=1,\ldots,N$. Then, for any $A \in \mathcal{F}$, we have
$$
P(A)=\sum_{j=1}^{N} P(A\,|\,B_j)P(B_j),
$$
where $P(A)$ is the total probability of event $A$.
\end{theorem}

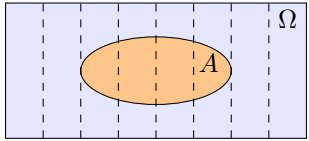
\begin{figure}[h] 
\centering 
\begin{tikzpicture}[xscale=2, yscale=0.9] 
\draw[fill=blue!10] (-1,-1)--(1,-1)--(1,1)--(-1,1)--(-1,-1);
\draw(0.875,0.5) node[above]{$\Omega$};
\draw[fill=orange!45] (0,0) circle (0.5);
\draw(0.35,-0.155) node[above]{$A$};
\draw[dashed] (-0.75,-1)--(-0.75,1);
\draw[dashed] (-0.5,-1)--(-0.5,1);
\draw[dashed] (-0.25,-1)--(-0.25,1);
\draw[dashed] (0,-1)--(0,1);
\draw[dashed] (0.25,-1)--(0.25,1);
\draw[dashed] (0.5,-1)--(0.5,1);
\draw[dashed] (0.75,-1)--(0.75,1);
\end{tikzpicture} 
\caption{Law of total probability} 
\label{Fig1}
\end{figure}

The main concept of Theorem \ref{TotalProbability} is the conditional probability. We define and denote the conditional probability of event $A$ given another event $B$ with $P(B)>0$ by
$$
P (A\,|\,B): =\frac{P(A \cap B)}{P(B)}.
$$

Conditional probability is a concept that allows us to interact with the reality of random phenomena by incorporating additional information about them.

On the other hand, Theorem \ref{TotalProbability} is crucial to understanding how probabilities propagate under different conditions, and its proof is straightforward as it relies on the basic properties of the probability measure. Indeed: Let $A \in \mathcal{F}$. We have,
$$
A=A \cap \Omega = A \cap \bigcup_{j=1}^{N} B_j = \bigcup_{j=1}^{N}(A \cap B_j)
$$
where $(A \cap B_{i}) \cap (A \cap B_{j})=\varnothing$ for all $i \neq j$ (see Figure \ref{Fig1}). Therefore, by the additivity of the probability measure, we obtain
$$
P(A)=\sum_{j=1}^{N} P(A \cap B_j) = \sum_{j=1}^{N} P(A\,|\,B_j)P(B_j).
$$

However, its essence is part of a more general concept in theory: \textit{conditional expectation} (see \cite[Definition 23.4]{Jacod(2004)}).

\begin{definition}\label{ConditionalExpectation}
Let $(\Omega,\mathcal{F},P)$ be a probability space, $X:\Omega \to \mathbb{R}$ (Lebesgue) integrable random variable, and $\mathcal{G}$ $\sigma$-algebra such that $\mathcal{G}\subset \mathcal{F}$. The conditional expectation of $X$ given $\mathcal{G}$ is the unique function $E(X\,|\,\mathcal{G}):\Omega \to \mathbb{R}$ with the following properties:
    \begin{itemize}
        \item[(i)] $E(X\,|\,\mathcal{G})$ is $\mathcal{G}$-measurable.
        \item[(ii)] $E(X\,|\,\mathcal{G})$ is (Lebesgue) integrable.
        \item[(iii)] For all $B \in\mathcal{G}$,
        $$
        \int_{B} X\,dP = \int_{B} E(X\,|\,\mathcal{G})\,dP.
        $$
    \end{itemize}
\end{definition}

Conditional expectation is fundamental for understanding the average behavior of a random variable, given certain conditions or additional information. This concept is closely related to the elementary ideas of probability. For the simplest case, when the $\sigma$-algebra consists only of the empty set and the entire sample space, $E(X\,|\,\mathcal{G})$ becomes a constant function equal to the unconditional expectation $E(X\,|\,\mathcal{G})=E(X)$. Consequently, if $X$ is an indicator function $1_{A}$ for an event $A \in \mathcal{F}$, then $E(1_{A}\,|\,\mathcal{G})$ equals the probability of $A$, $P(A)$. 

The above suggests that the basic concept of conditional probability can be derived from conditional expectation. To do this, we need to carefully study a certain class of $\sigma$-algebras.

The idea to follow in the text is as follows: Let $(\Omega,\mathcal{F},P)$ be a probability space and $\mathcal{B}=\{B_1,\ldots,B_N\}$  a collection of elements from $\mathcal{F}$ such that they form a partition of $\Omega$ and $P(B_j)>0$ for all $j=1,\ldots,N$. 

We define $\sigma$-algebra $\mathcal{G}$ as the $\sigma$-algebra generated by the class $\mathcal{B}$, that is, $\mathcal{G}:=\sigma(\mathcal{B})$. It is clear that $\mathcal{G} \subset \mathcal{F}$.

Fixed $A \in \mathcal{F}$. We define the measure $Q:\mathcal{G} \to \mathbb{R}$
$$
Q(B):=\int_{B} 1_{A}\,dP
$$
where it is clear that $Q(B) \leq P(B)$ for all $B \in \mathcal{G}$ and therefore $Q$ is a finite measure. Our goal is to prove the existence of a unique $\mathcal{G}$-measurable and Lebesgue integrable function $\xi:\Omega \to \mathbb{R}$ such that
\begin{equation} \label{F1}
    Q(B)=\int_{\Omega} \xi\cdot 1_{B}\,dP\quad \mbox{for all}\,\,B \in \mathcal{G}.
\end{equation}

That is, we want to prove the existence of the conditional expectation of $1_{A}$ given $\mathcal{G}$ and deduce the law of total probability.

We begin establishing notation and useful facts in Section \ref{usefulfacts}, where we provide a brief overview of different results from functional analysis and probability, as we use them interchangeably in the developed sections. These results include original proofs that are tailored to the intended objective.

In Section \ref{Statement}, we focus on providing an overview of the problem to be solved, that is, we will frame the existence and uniqueness of the conditional expectation of an event A as a linear problem in a vector space of infinite dimension.

In Section \ref{L2}, we restrict the problem to the Hilbert space $L^{2}$ of square-integrable random variables, where the inner product is given by the expectation of the product of two random variables. This restriction allows us to use the Riesz-Fréchet representation theorem and, therefore, to provide a geometric description of the conditional expectation. In the literature (see \cite[Chapter 23]{Jacod(2004)} for example), it is common to prove the existence of the conditional expectation in $L^{2}$ with this description. However, in this work, we point out that this is not intrinsic to the $L^2$ structure. In fact, in Section \ref{Lp} we show that for certain values of $p \in (1,\infty)$ the geometry of the space $L^p$ of random variables with the $p$-th power absolutely integrable is similar to that of $L^2$. We do this by presenting a version of the Riesz representation theorem and using Clarkson's inequalities.

In Section \ref{VariationalProblem}, we frame the existence of the conditional expectation as the existence of a critical point of a certain functional defined in an infinite-dimensional Hilbert space. This approach is known as the Dirichlet principle and its solution is given by the Riesz-Fréchet representation theorem. Following this idea and through the Hölder-Riesz inequality, in this paper we extend this principle to certain values of $p \in (1,\infty)$.

Finally, in Section \ref{TP}, where we prove the law of total probability (Theorem \ref{TotalProbability}) using the tools developed in the previous sections. This proof, unlike the one given in the initial paragraphs, relies only on properties of the Lebesgue integral.

\begin{remark} In this text, we consider two random variables identical if they are equal except on a set of probability zero. That is, $X=Y$ if $P(X=Y)=P(\{\omega \in \Omega \,:\, X(\omega)=Y(\Omega)\})=1$.
\end{remark}

\section{Notation and some useful facts} \label{usefulfacts}

The main idea of the text is to frame (\ref{F1}) as a linear problem in certain Banach spaces. Therefore, in this section, we will detail some key concepts and results from functional analysis and probability theory. We will also establish some notation.

\subsection[\appendixname~\thesubsection]{Banach spaces}

Let $V=(V,\Vert \cdot \Vert_V)$ and $W=(W,\Vert \cdot \Vert_W)$ be two Banach spaces over $\mathbb{R}$.

\begin{definition}\label{Continuity}
    A function $f:V \to W$ is continuous on $V$ if, for every $u \in V$ and $(u_j)$ sequence on $V$ such that $\lim_{j \to \infty}\Vert u_j - u \Vert_{V}=0$, then $\lim_{j \to \infty} \Vert f(u_j) - f(u) \Vert_{W}=0$.
\end{definition}

\begin{definition} \label{CFunction} For linear and bilinear functions, the following characterization of continuity holds:

    \begin{itemize}
    \item[(a)] A linear function $L:V \to W$ is continuous if exists a constant $C >0$ such that $\Vert L(v) \Vert_{W} \leq C \Vert v \Vert_{V}$ for all $v \in V$.

    \item[(b)] A bilinear function $L:V \times V \to W$ is continuous if there exists a constant $M >0$ such that $\Vert L(u,v) \Vert_{W} \leq M \Vert u \Vert_{V}\Vert v \Vert_{V}$ for all $u,v \in V$.
    \end{itemize}
\end{definition}

The space:
$$
\mathcal{B}(V,W):=\left\{L:V \to W\,:\,L\,\,\mbox{is linear and continuos} \right\}
$$
with norm:
$$
\Vert L \Vert_{\mathcal{B}(V,W)}:=\sup_{\Vert v \Vert_{V}\leq 1} \Vert L(v) \Vert_{W}
$$
is a Banach space over $\mathbb{R}$ \cite[Theorem 2.10-2]{Krey(1978)}. If $W=\mathbb{R}$ we denote this space like $V^{\ast}:=\mathcal{B}(V,\mathbb{R})$ and we call it the dual space of $V$. In the same way, we define the space $V^{\ast\ast}:=\mathcal{B}(V^{\ast},\mathbb{R})$ and call it the bidual space of $V$.

If $L \in \mathcal{B}(V,W)$, then $\Vert L(v) \Vert_{W} \leq \Vert L \Vert_{\mathcal{B}(V,W)}$ for all $\Vert v \Vert_{V} \leq 1$. Consequently, $\Vert L(v) \Vert_{V} \leq \Vert L \Vert_{\mathcal{B}(V,W)} \Vert v \Vert_{V}$ for all $v \in V$.

Let $v \in V$. We define function $\widehat{v}:V^{\ast} \to \mathbb{R}$ by
$$
\widehat{v}(L):=L(v).
$$

It is clear that $\hat{v}$ is linear in $V^{\ast}$. Note that $|L(v)| \leq \Vert L \Vert_{V^{\ast}} \Vert v \Vert_{V}$ for all $L \in V^{\ast}$ and therefore $|\widehat{v}(L)| = |L(v)| \leq \Vert L \Vert_{V^{\ast}} \Vert v \Vert_{V}$ for all $L \in V^{\ast}$. This implies that $\widehat{v}$ is continuous at $V^{\ast}$.

The following result establishes a relationship between the space $V$ and its bidual space \cite[Lemma 4.6.2]{Krey(1978)}.

\begin{theorem}\label{CanonicalIsometry}
The function $\mathcal{I}:V \to V^{\ast \ast}$ defined by
$$
\mathcal{I}(v):=\widehat{v}
$$
is linear, continuous, and $\Vert \mathcal{I}(v) \Vert_{V^{\ast \ast}}=\Vert v \Vert_{V}$ for all $v\in V$.
\end{theorem}

The function $\mathcal{I}$ in Theorem \ref{CanonicalIsometry} is called the canonical mapping of $V$ into $V^{\ast \ast}$ and we can deduce that it is injective.

\begin{definition}\label{ReflexiveSpaces}
    A Banach space $V$ is said to be reflexive if the canonical mapping $\mathcal{I}$ is surjective.
\end{definition}

The following concept will provide us with a sufficient condition for when a space is reflexive, the Milman-Pettis theorem \cite[Theorem 3.31]{Brezis(2011)}.

\begin{definition} \label{UniformConvexity}
A normed space $V=(V,\Vert \cdot \Vert_{V})$ is said to be uniformly convex if, given an arbitrary real number $\varepsilon >0$, there is a $\delta>0$ such that for any $v,w \in V$ with $\Vert v \Vert_V \leq 1$, $\Vert  w \Vert_V \leq 1$ and $\Vert v-w\Vert_V > \varepsilon$ we have
$$
\left\Vert \frac{v+w}{2} \right\Vert_{V} < 1-\delta.
$$
\end{definition}

Essentially, a normed space is uniformly convex if its closed unit ball does not contain line segments on its boundary. The most precise example of this is a Hilbert space.

\begin{theorem}[Milman-Pettis] \label{MilmanPettis}
Let $V$ be a uniformly convex Banach space. Then $V$ is reflexive.
\end{theorem}

\begin{corollary} \label{ReflexiveHilbert}
    If $H=(H,\langle \cdot, \cdot \rangle, \Vert \cdot \Vert_{H})$ is the Hilbert space over $\mathbb{R}$, then $H$ is reflexive.
\end{corollary}

\begin{proof}
    Let $\varepsilon >0$ and $u,v \in H$ be such that $\Vert u \Vert_{H},\Vert v \Vert_{H} \leq 1$ and $\Vert u-v \Vert_{H} > \varepsilon$. We have,
    $$
    \begin{aligned}
    \Vert u+v\Vert_{H}^{2}&=\langle u+v,u+v\rangle = \Vert u \Vert_{H}^2 +2\langle u,v\rangle + \Vert v \Vert_{H}^2,\\
    \Vert u-v\Vert_{H}^{2}&=\langle u-v,u-v\rangle = \Vert u \Vert_{H}^2 -2\langle u,v\rangle + \Vert v \Vert_{H}^2,
    \end{aligned}
    $$
    and thus $\Vert u+v\Vert_{H}^{2} + \Vert u-v\Vert_{H}^{2}=2(\Vert u \Vert_{H}^{2} + \Vert v \Vert_{H}^{2})$. Hence,
    $$
    \left\Vert \frac{u+v}{2} \right\Vert_{H}^{2} = \frac{1}{2}\left( \Vert u \Vert_{H}^{2} + \Vert v \Vert_{H}^{2} \right)-\left\Vert \frac{u-v}{2} \right\Vert_{H}^{2} < 1 - \frac{\varepsilon^2}{4}
    $$
    with $\delta:=1-(1-\frac{\varepsilon^2}{4})^{1/2} >0$ we have the result, by Theorem \ref{MilmanPettis}.
\end{proof}

\subsection{Probability spaces} 

A probability space $(\Omega,\mathcal{F},P)$ consists of three elements:
\begin{itemize}
    \item[(a)] $\Omega$ is a non-empty set.
    \item[(b)] $\mathcal{F}$ is a $\sigma$-algebra in $\Omega$, that is, $\mathcal{F}$ is a collection of subsets of $\Omega$ such that:
        \begin{itemize}
            \item[(b.1)] $\Omega \in \mathcal{F}$,
            \item[(b.2)] If $A \in \mathcal{F}$ then $\Omega \smallsetminus A \in \mathcal{F}$,
            \item[(b.3)] $\bigcup_{k=1}^{\infty} A_k \in \mathcal{F}$ whenever $A_k \in \mathcal{F}$ for all $k \in \mathbb{N}$. 
        \end{itemize}
    \item[(c)] $P$ is a probability measure, that is, $P:\mathcal{F} \to [0,1]$ satisfies:
    \begin{itemize}
            \item[(c.1)] $P(\Omega)=1$,
            \item[(c.2)] $P\left(\bigcup_{k=1}^{\infty} A_k\right)=\sum_{k=1}^{\infty} P(A_k)$ whenever $(A_k)$ is a disjoint countable family of members of $\mathcal{F}$. 
        \end{itemize}
\end{itemize}

\begin{definition} \label{SigmaGenerated}
    Let $\mathcal{C}$ be an arbitrary collection of subsets of de $\Omega$. We define the $\sigma$-algebra generated by $\mathcal{C}$, denoted by $\sigma(\mathcal{C})$, as the intersection of all the $\sigma$-algebras of subsets of $\Omega$ that contain $\mathcal{C}$, that is,
$$
\sigma(\mathcal{C}):=\bigcap\{\mathcal{F}\,:\,\mathcal{F}\,\,\sigma\,\mbox{-algebra and}\,\,\mathcal{C}\subset \mathcal{F}\}.
$$
\end{definition}

$\sigma(\mathcal{C})$ is a $\sigma$-algebra of subsets of $\Omega$ and, in fact, it is the smallest $\sigma$-algebra that contains $\mathcal{C}$. Indeed, by definition, it is clear that $\mathcal{C}\subset \sigma(\mathcal{C})$ and if $\mathcal{F}_0$ is a $\sigma$-algebra such that $\mathcal{C} \subset \mathcal{F}_0$ then $\mathcal{F}_0 \in \{\mathcal{F}\,:\,\mathcal{F}\,\,\sigma\,\mbox{-algebra and}\,\,\mathcal{C}\subset \mathcal{F}\}$ and, thus, $\sigma(\mathcal{C})\subset \mathcal{F}_0$.

It is possible to provide a description of the $\sigma$-algebra generated by a finite partition of $\Omega$.

\begin{proposition} \label{PropSigmaG}
    Let $\mathcal{B}=\{B_1,\ldots,B_N\}$ be a finite collection of elements from $\mathcal{F}$ such that they form a partition of $\Omega$. Then
$$
\sigma(\mathcal{B})=\left\{ \bigcup_{k \in \tau} B_k\,:\, \tau \subset \{1,\ldots,N\}\,\,\mbox{and}\,\,\tau \neq \varnothing \right\}\cup \{\varnothing\}.
$$
\end{proposition}

\begin{figure}[h!]
\centering 
\begin{tikzpicture}[xscale=1.75, yscale=1.2] 

\draw[fill=blue!10] (-1,-1)--(1,-1)--(1,1)--(-1,1)--(-1,-1);
\draw(0.875,0.62) node[above]{$\Omega$};

\draw[dashed] (-0.75,-1)--(-0.75,1);
\draw[dashed] (-0.5,-1)--(-0.5,1);
\draw[dashed] (-0.25,-1)--(-0.25,1);
\draw[dashed] (0,-1)--(0,1);
\draw[dashed] (0.25,-1)--(0.25,1);
\draw[dashed] (0.5,-1)--(0.5,1);
\draw[dashed] (0.75,-1)--(0.75,1);
\end{tikzpicture} 
\caption{Partition of $\Omega$}\label{Fig2}
\end{figure}
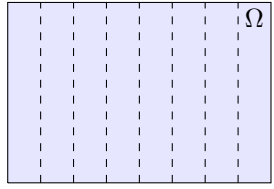

\begin{proof}
    We define
$$\mathfrak{A}:=\left\{ \bigcup_{t \in \tau} B_t\,:\, \tau \subset \{1,\ldots,N\}\,\,\mbox{and}\,\,\tau \neq \varnothing \right\}\cup \{\varnothing\}.$$ 

It is clear that $\mathcal{B} \subset \mathfrak{A}$ and by definition $\sigma$-algebra, $\mathfrak{A} \subset \sigma(\mathcal{B})$. We will prove that $\mathfrak{A}$ is a $\sigma$-algebra. Note that since the elements of $\mathcal{B}$ form a finite partition of $\Omega$, $\mathfrak{A}$ is a finite set such that it suffices to prove that $\mathfrak{A}$ is closed under finite unions.

Given by $\Omega=B_1\cup\cdots\cup B_N$ then $\Omega \in \mathfrak{A}$. Now, let $A \in \mathfrak{A}$. It is clear that if $A=\varnothing$ then $\Omega\smallsetminus A = \Omega \in \mathfrak{A}$ and in the same way if $A=\Omega$ then $\Omega \smallsetminus A = \varnothing \in \mathfrak{A}$. We suppose that $A \neq \varnothing$ and that there exists a proper non-empty subset $\tau$ of $\{1,\ldots,N\}$ such that
$$
A=\bigcup_{t \in \tau} B_{t}.
$$

Again, since $\Omega = \bigcup_{j \in \{1,\ldots,N\}} B_{j}$ then
$$
\Omega \smallsetminus A = \left(\bigcup_{j \in \{1,\ldots,N\}} B_{j} \right) \smallsetminus \left(\bigcup_{t \in \tau} B_{t} \right) = \bigcup_{t \in \{1,\ldots,N\} \smallsetminus \tau } B_{t} \in \mathfrak{A}.
$$

Finally, if $A_1,A_2 \in \mathfrak{A}\smallsetminus \{\varnothing\}$, there exist nonempty subsets $\tau_1$ and $\tau_2$ of \\$\{1,\ldots,N\}$ such that
$$
A_1=\bigcup_{i \in \tau_1} B_{i}\,\,\,\,\,\,\,\,\mbox{and}\,\,\,\,\,\,\,\,A_2=\bigcup_{j \in \tau_2} B_{j}.
$$

Therefore,
$$
A_1 \cup A_2 = \left( \bigcup_{i \in \tau_1} B_{i} \right) \cup \left(\bigcup_{j \in \tau_j} B_{j}\right) = \bigcup_{t \in \tau_1 \cup \tau _2} B_{t},
$$
and hence $A_1 \cup A_2 \in \mathfrak{A}$. The case where $A_{1}=\varnothing$ or $A_{2}=\varnothing$ is clear. Consequently, $\mathfrak{A}$ is a $\sigma$-algebra in $\Omega$ such that $\mathcal{B} \subset \mathfrak{A}$ and thus $\sigma(\mathcal{B})\subset \mathfrak{A}$. 
\end{proof}

\begin{definition} \label{RandomVariable}
    $X:\Omega \to \mathbb{R}$ is a random variable (or a $\mathcal{F}$-measurable function) if \\ $X^{-1}(-\infty,c] \in \mathcal{F}$ for all $c \in \mathbb{R}$\end{definition}

In the previous definition, it is important to note that the measurability of the function depends on the $\sigma$-algebra considered in the probability space.

\begin{proposition} \label{IndicatorFunction}
    If $A \in \mathcal{F}$ then indicator function $1_{A}$ given by 
    $$
    1_{A}(\omega):=\left\{
    \begin{array}{c l}
    1 & if\,\,\omega \in A,\\
    0 & if\,\, \omega\not\in A,
    \end{array}
\right.
$$
is a random variable.
\end{proposition}

\begin{proof}
    For $c \in \mathbb{R}$ we have that
    $$
1_{A}^{-1}(-\infty,c]=\{\omega \in \Omega \,:\, 1_{A}(\omega) \leq c \}=\left\{
\begin{array}{l c l}
\Omega & & \mbox{if}\,\, c \geq 1,\\
\Omega \smallsetminus A & & \mbox{if}\,\, 0 \leq c <1,\\
\varnothing & & \mbox{if}\,\, c<0.
\end{array}
\right.
$$

Thus, $1_{A}^{-1}(-\infty,c] \in \mathcal{F}$ for all $c \in \mathbb{R}$.
\end{proof}

\begin{proposition} \label{IndicatorLinealCombination}
    A linear combination of indicator functions of sets in $\mathcal{F}$ is a random variable.
\end{proposition}

\begin{figure}[h!]
    \centering
	\begin{tikzpicture}[xscale=0.8,yscale=0.8]
	\draw[->,gray] (0,0) -- (4,0); \draw [->,gray] (0,0) -- (0,4); 
	 \draw[-,gray] (0,0)--(0,-1.2);

\draw (0,2.8) node{$-$}; \draw (0,2.8) node[left] {$_{\alpha_{i}}$};
\draw [ thick] (1.3,2.8)--(2.22,2.8);
\draw (1.3,2.8) node{$\bullet$}; 
\draw (2.3,2.8) node{$\circ$};

\draw (0,2.3) node{$-$}; \draw (0,2.3) node[left] {$_{\alpha_{k}}$};
\draw [ thick] (2.3,2.3)--(2.9,2.3);
\draw (2.3,2.3) node{$\bullet$}; 
\draw (3,2.3) node{$\circ$};

\draw (0,1) node{$-$}; \draw (0,1) node[left] {$_{\alpha_{j}}$};
\draw [ thick] (0,1)--(1.12,1);
\draw (0,1) node{$\bullet$}; 
\draw (1.2,1) node{$\circ$};

\draw (0,-1) node{$-$}; \draw (0,-1) node[left] {$_{\alpha_{\ell}}$};
\draw [ thick] (3,-1)--(4,-1);
\draw (3,-1) node{$\bullet$}; 

\draw[densely dotted, gray] (1.25,2.8)--(1.25,1.1);
\draw[densely dotted, gray] (1.25,0.9)--(1.25,0);

\draw[densely dotted, gray] (2.3,2.7)--(2.3,0); 
\draw[densely dotted, gray] (3,-1)--(3,2.25); 
\end{tikzpicture}
\caption{Linear combination of indicator functions}
\label{Fig3}
\end{figure}
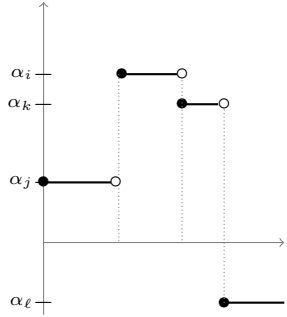

\begin{proof}
    Let $X=\sum_{j=1}^{N} \alpha_{j}1_{A_{j}}$ with $\alpha_{j} \in \mathbb{R}$ and $A_{j} \in \mathcal{F}$ for all $j =1,\ldots,N$ (see Figure \ref{Fig3}).

    If $c \in \mathbb{R}$ is any real number, then
    $$
     \{\omega \in \Omega\,:\, X(\omega) \leq  c \} = \left\{
    \begin{array}{lcl}
    \varnothing & & \mbox{if}\,\,\alpha_j  > c\,\,\forall\,\,j=1,\ldots,N,\\
    \bigcup\{A_j\,:\, \alpha_j \leq c \} & & \mbox{another case}.
    \end{array}
    \right.
    $$

    Consequently, $X^{-1}(-\infty,c] \in \mathcal{F}$ for all $c \in \mathbb{R}$.
\end{proof}

A function that is a linear combination of indicator functions is referred to as a simple random variable \cite[Definition 9.1]{Jacod(2004)}.

With the previous results, we can characterize the measurable functions in the $\sigma$ algebra generated by a finite partition of $\Omega$.

\begin{proposition} \label{GMeasurable}
Let $\mathcal{B}=\{B_1,\ldots,B_N\}$ be a finite collection of elements from $\mathcal{F}$ such that they form a partition of $\Omega$. If $\mathcal{G}:=\sigma(\mathcal{B})$ then $X:\Omega \to \mathbb{R}$ is $\mathcal{G}$ measurable if and only if for each $j=1,\ldots,N$ there exists $\alpha_{j} \in \mathbb{R}$ such that $X=\sum_{j=1}^{N} \alpha_{j}1_{B_j}$.
\end{proposition}

\begin{proof}
$\Leftarrow):$ It's the Proposition \ref{IndicatorLinealCombination}.

$\Rightarrow):$ Suppose, by contradiction, that there exist $j_0 \in \{1,\ldots,N\}$ and $\omega_{\alpha},\omega_{\beta} \in B_{j_{0}}$ such that $X(\omega_{\alpha})=\alpha < \beta = X(\omega_{\beta})$.

We know that $X$ is $\mathcal{G}$-measurable, so $X^{-1}(-\infty,\alpha] \in \mathcal{G}$ and therefore there exists $\tau \subset \{1,\ldots,N \}$ such that $X^{-1}(-\infty,\alpha]=\bigcup_{t \in \tau} B_{t}$ by Proposition \ref{PropSigmaG}. It is clear that $B_{j_{0}} \cap X^{-1}(-\infty,\alpha] \neq \varnothing$ and using the hypothesis, we obtain $B_{j_{0}} \cap X^{-1}(-\infty,\alpha] =\bigcup_{t \in \tau} (B_{t} \cap B_{j_{0}})$. However, the sets $B_{j}$ form a partition of $\Omega$ so $B_{j_{0}} \cap X^{-1}(-\infty,\alpha] \neq \varnothing$ if and only if $B_{j_{0}} \cap X^{-1}(-\infty,\alpha]= B_{j_{0}}$. This is a contradiction, since $\omega_{\beta} \in B_{j_{0}} \smallsetminus X^{-1}(-\infty,\alpha]$.
\end{proof}

\begin{definition} \label{Expectation}
    If $X:\Omega \to \mathbb{R}$ is a random variable in $(\Omega,\mathcal{F},P)$ we define the expectation using the Lebesgue integral, that is,
    $$
    E(X):=\int_{\Omega} X\,dP.
    $$

    We say that $X$ is (Lebesgue) integrable if $E(|X|)<\infty$.
\end{definition}

In this paper, we consider two random variables identical if they are equal except on a set of probability zero. That is, $X=Y$ if $P(X=Y)=P(\{\omega \in \Omega \,:\, X(\omega)=Y(\Omega)\})=1$.

We write $L^{1}(\Omega,\mathcal{F},P)$ to denote the set of all integrable random variables on $(\Omega,\mathcal{F},P)$. We state some properties of expectation whose proofs can be found in \cite[Theorem 9.1]{Jacod(2004)}.

\begin{theorem} \label{L1vectorspace}
Let $(\Omega,\mathcal{F},P)$ be a probability space.
\begin{itemize}
    \item[(a)] $L^{1}(\Omega,\mathcal{F},P)$ is a vector space over $\mathbb{R}$, and expectation is a linear map on $L^{1}(\Omega,\mathcal{F},P)$, that is, if $X,Y \in L^{1}(\Omega,\mathcal{F},P)$ and $\alpha,\beta \in \mathbb{R}$ then $E(\alpha X+ \beta Y)=\alpha E(X)+\beta E(Y)$.
    \item[(b)] $X \in L^{1}(\Omega,\mathcal{F},P)$ if and only if $|X| \in L^{1}(\Omega,\mathcal{F},P)$ and in this case $|E(X)| \leq E(|X|)$.
    \item[(c)] If $X,Y \in L^{1}(\Omega,\mathcal{F},P)$ then $X=Y$ if and only if $E(X)=E(Y)$. 
    \item[(d)] If $X,Y \in L^{1}(\Omega,\mathcal{F},P)$ is such that $X \leq Y$ then $E(X) \leq E(Y)$
\end{itemize}
\end{theorem}

For $X \in L^{1}(\Omega,\mathcal{F},P)$ we denote:
$$
\Vert X \Vert_{1}:=E(|X|).
$$

Thus, $L^{1}(\Omega,\mathcal{F},P)$ with the norm $\Vert \cdot \Vert_{1}$ is a Banach space. This result is known as the Riesz-Fischer theorem, and its proof is provided in \cite[Theorem 4.8]{Brezis(2011)}.

One can show that if $X$ is a random variable then there exists a sequence $(X_k)$ of simple random variables such that $|X_k| \leq |X|$ for any $k \in \mathbb{N}$ and $X_k \to X$ on $\Omega$ (see \cite[Result 1, Chapter 9]{Jacod(2004)}). Hence, if $X \in L^1(\Omega,\mathcal{F},P)$ then there exists a sequence $(X_k)$ of simple random variables in $L^{1}(\Omega,\mathcal{F},P)$ such that $\lim_{k \to 0}\Vert X_k - X\Vert_{1}=0$ (see \cite[Theorem 9.1]{Jacod(2004)}).

\section{Problem statement} \label{Statement}

Through the results of the previous results in Section \ref{usefulfacts}, we present the main idea of the problem addressed in this paper.

Let $(\Omega,\mathcal{F},P)$ be a probability space and $\mathcal{B}=\{B_1,\ldots,B_N\}$ a finite collection of elements from $\mathcal{F}$ that form a partition of $\Omega$. Following the idea of the law of total probability, we interpret the collection $\mathcal{B}$ as additional information on our random phenomenon. 

We consider $\sigma$-algebra $\mathcal{G}:=\sigma(\mathcal{B})$ and it is easy to see that $\mathcal{G} \subset \mathcal{F}$. If we take the restriction of $P$ to $\mathcal{G}$, we can conclude that $(\Omega,\mathcal{G},P)$ is a probability space.

Fixed $A \in \mathcal{F}$, our problem is to find a unique function $\xi \in L^{1}(\Omega,\mathcal{G},P)$ such that:
\begin{equation} 
    \int_{\Omega} 1_{A}\cdot 1_{B}\,dP =\int_{\Omega} \xi\cdot 1_{B}\,dP\,\,\,\,\,\,\mbox{for all}\,\,B \in \mathcal{G}.
\end{equation}

Note that $1_{B} \in L^{1}(\Omega,\mathcal{G},P)$ for all $B \in \mathcal{G}$ so that we can extend the idea as follows: 

We define $T:L^{1}(\Omega,\mathcal{G},P) \to\mathbb{R}$ by:
$$
T(X):=\int_{\Omega} X\cdot 1_{A}\,dP.
$$

The items $(a)$ and $(b)$ of Theorem \ref{L1vectorspace} ensure that $T$ is a continuous linear function on $L^{1}(\Omega,\mathcal{G},P)$. Then, our problem reduces to find a unique function $\xi \in L^{1}(\Omega,\mathcal{G,P})$ such that
\begin{equation} 
    T(X) =\int_{\Omega} \xi\cdot X\,dP\,\,\,\,\,\,\mbox{for all}\,\,X \in L^{1}(\Omega,\mathcal{G},P).
\end{equation}

The above is similar to a Riesz representation theorem, so in the following sections, we will study it in detail.

\section{Problem in $L^2$} \label{L2}
We consider $(\Omega,\mathcal{G},P)$ the probability space given in Section \ref{usefulfacts} and write
$$
L^{2}(\Omega,\mathcal{G},P):=\left\{ X:\Omega \to \mathbb{R}\,:\,X\,\,\mbox{is a random variable and}\,\,E(|X|^2)< \infty \right\}.
$$

It is easy to see $X \in L^{2}(\Omega,\mathcal{G},P)$ if and only if $|X|^{2} \in L^{1}(\Omega,\mathcal{G},P)$. Hence, $L^{2}(\Omega,\mathcal{G},P)$ is a vector space over $\mathbb{R}$.

The fundamental result in this section is the H\"older-Riesz inequality \cite[Theorem 4.6]{Brezis(2011)}, as it will allow us to tackle the initial problem in this space.

Previously, for $X \in L^{2}(\Omega,\mathcal{F},P)$ we denote
\begin{equation}
    \Vert X \Vert_{2}:=\left( E(|X|^2) \right)^{1/2}.
\end{equation}

\begin{theorem}[H\"older-Riesz inequality]  \label{HR2}
If $X, Y \in L^{2}(\Omega,\mathcal{G},P)$, then $X\cdot Y \in L^{1}(\Omega)$ and $\Vert X \cdot Y \Vert_{1} \leq \Vert X \Vert_{2} \Vert Y \Vert_{2}$.
\end{theorem}

We will also refer to the following corollary as the H\"older-Riesz inequality.

\begin{corollary} \label{CorHR2}
If $X \in L^{2}(\Omega,\mathcal{G},P)$, then $X \in L^{1}(\Omega,\mathcal{G},P)$ and $\Vert X \Vert_{1} \leq \Vert X \Vert_{2}$.
\end{corollary}

\begin{proof}
    Apply Theorem \ref{HR2} to the functions $X$ and $1_{\Omega}$.
\end{proof}

In $L^{2}(\Omega,\mathcal{G},P)$ we define
$$
\langle X,Y\rangle_{2}:=E(X\cdot Y).
$$

By the H\"older-Riesz inequality and Theorem \ref{L1vectorspace}, we claim that it is an inner product in $L^{2}(\Omega,\mathcal{G},P)$ such that
$$
\langle X,X \rangle_{2}^{2}=\Vert X \Vert_{2}^{2}.
$$

In fact, $L^{2}(\Omega,\mathcal{G},P)$ is a Hilbert space because $L^{2}(\Omega,\mathcal{G},P)$ with norm $\Vert \cdot \Vert_{2}$ is a Banach space \cite[Theorem 4.8]{Brezis(2011)}.

We rewrite the initial problem in $L^{2}(\Omega,\mathcal{G},P)$ as follows: 
\\Define $T:L^{2}(\Omega,\mathcal{G},P) \to \mathbb{R}$ by
\begin{equation} \label{FunctionalT2}
    T(X):=\int_{\Omega} X\cdot 1_{A}\,dP.
\end{equation}

Again, by Theorem \ref{L1vectorspace} $T$ is linear in $L^{2}(\Omega,\mathcal{G},P)$ and by the H\"older-Riesz inequality, $|T(X)| \leq \Vert X \Vert_{2}$ for all $X \in L^{2}(\Omega,\mathcal{G},P)$, which implies that $T$ is continuous in $L^{2}(\Omega,\mathcal{G},P)$ (see Definition \ref{CFunction}). The question is: Is there a unique function $\xi \in L^{2}(\Omega,\mathcal{G},P)$ such that $T(X)=E(\xi \cdot X)$ for all $X \in L^{2}(\Omega,\mathcal{G},P)$?

In Hilbert spaces, there is a result that allows us to represent the elements of their dual through their inner product, the Riesz–Fréchet representation theorem \cite[Theorem 3.8-1]{Krey(1978)}.

\begin{theorem}[Riesz–Fréchet representation theorem] \label{FrechetRiesz}
Let $H=(H,\langle\cdot,\cdot \rangle,\Vert \cdot \Vert_{H})$ be a Hilbert space over $\mathbb{R}$ and $L:H \to \mathbb{R}$ a continuous linear function. There exists a unique $u \in H$ such that:
$$
L(v)=\langle u,v\rangle \,\,\,\,\,\,\,\forall\,v \in H.
$$
\end{theorem}

Hence, by Theorem \ref{FrechetRiesz} there exists a unique $\xi \in L^{2}(\Omega,\mathcal{G},P)$ such that
$$
T(X)=\langle \xi,X \rangle_{2} = \int_{\Omega} \xi\cdot X\,dP \,\,\,\,\,\,\forall\,\,X \in L^{2}(\Omega,\mathcal{G},P).
$$

Now, by Corollary \ref{HR2} $\xi \in L^{1}(\Omega,\mathcal{G},P)$ and we have that
$$
T(1_B)=\int_{\Omega} \xi\cdot 1_{B}\,dP \,\,\,\,\,\,\forall\,\,B \in \mathcal{G}
$$
because $1_{B} \in L^{2}(\Omega,\mathcal{G},P)$ for all $B \in \mathcal{G}$. Consequently, $\xi$ is $E(1_{A}\,|\,\mathcal{G})$.

The question now is, how do we interpret $E(1_{A}\,|\,\mathcal{G})$ through the proof of Theorem \ref{FrechetRiesz}?

\begin{figure}[h!] 
\centering 
\begin{tikzpicture}[xscale=1.55, yscale=1.55] 

\draw[->, thick] (-1,0)--(1,0); 
\draw[->, thick] (0,-1)--(0,1); 

\draw(1.35,0) node{$\mbox{ker\,L}$};
\draw(0,1.1) node{$(\mbox{ker\,L})^{\perp}$};

\draw(0,0.35) node[right]{$_{u}$}; 

\draw[->, blue] (0,0)--(0,0.35);
\end{tikzpicture} 
\caption{Orthogonal decomposition of $H$ in the kernel of $L$}
\label{Fig4}
\end{figure}
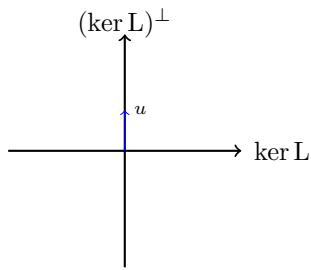

In \cite{Krey(1978)}, the idea of the proof of Theorem \ref{FrechetRiesz} is: Since $L:H \to \mathbb{R}$ is a linear and continuous function, then $\mbox{ker}L=L^{-1}(\{0\})$ is a closed vector subspace of $H$ and, thus, $H=\mbox{ker}L \oplus (\mbox{ker}L)^{\perp}$ (see Figure \ref{Fig4}).

If $L$ is trivial, then $\mbox{ker}L=H$ and $u=0_{H}$ is a vector zero. If $L$ is nontrivial, then $(\mbox{ker}\,L)^{\perp} \neq \{0\}$. Fix $w \in \mbox{ker}L$ with $\Vert w \Vert_H=1$ and defining $u:=L(w)w \in H$ yields the result.

That is, we can interpret $E(1_{A}\,|\,\mathcal{G})$ as an element of $(\mbox{ker}\,T)^{\perp}$.

However, in the proof of Theorem \ref{FrechetRiesz} given in \cite[Theorem 5.5]{Brezis(2011)}, two important tools are used: the Hahn-Banach theorem \cite[Corollary 1.8]{Brezis(2011)} and the reflexivity of Hilbert spaces (see Corollary \ref{ReflexiveHilbert}).

We then ask: Is orthogonality necessary for our problem? Is there a similar theorem in spaces that are not necessarily Hilbert spaces? In the following section, we will answer these questions.

\section{Problem in $L^{p}$} \label{Lp}

We find that there are various Riesz representation theorems, and in this section we will focus on $L^p$ spaces. What is special about Lebesgue spaces $L^p$ is that the geometric properties we are interested in are similar to those in the Hilbert space $L^2$. And this is enough to tackle the problem.

As before, we consider $(\Omega,\mathcal{G},P)$ the probability space given in Section \ref{usefulfacts} and write, for $p \in (1,\infty)$
$$
L^{p}(\Omega,\mathcal{G},P):=\left\{ X:\Omega \to \mathbb{R}\,:\,X\,\,\mbox{is a random variable and}\,\,E(|X|^{p})< \infty \right\}.
$$

It's clear that $X \in L^{p}(\Omega,\mathcal{G},P)$ if and only if $|X|^{p} \in L^{1}(\Omega,\mathcal{G},P)$. Hence, $L^{p}(\Omega,\mathcal{G},P)$ is a vector space over $\mathbb{R}$ \cite[Theorem 4.7]{Brezis(2011)}.

For $X \in L^{p}(\Omega,\mathcal{F},P)$ we denote
\begin{equation}
    \Vert X \Vert_{p}:=\left( E(|X|^p) \right)^{1/p}.
\end{equation}

Let $p \in (1,\infty)$. We denote by $q \in (1,\infty)$ the conjugate exponent, that is, $q:=\frac{p}{p-1}$ and it is such that
$$
\frac{1}{p}+\frac{1}{q}=1.
$$

There is a general version of the H\"older-Riesz inequality \cite[Theorem 4.6]{Brezis(2011)}, which we state as follows.

\begin{theorem}[H\"older-Riesz inequality in $L^{p}$] \label{HRp}
If $p,q\in (1,\infty)$ are conjugate exponents, $X\in L^{p}(\Omega,\mathcal{G},P)$ and $Y\in L^{q}(\Omega,\mathcal{G},P)$ then $X\cdot Y \in L^{1}(\Omega)$ and $\Vert X \cdot Y \Vert_{1} \leq \Vert X \Vert_{p} \Vert Y \Vert_{q}$.
\end{theorem}

\begin{corollary} \label{CHRp}
    Let $1 \leq r < s < \infty$. Then $L^{s}(\Omega,\mathcal{G},P) \subset L^{r}(\Omega,\mathcal{G},P)$. In addition, for any $X \in L^{s}(\Omega,\mathcal{G},P)$ we have $\Vert X \Vert_{r} \leq \Vert X \Vert_{s}$
\end{corollary}

\begin{proof}
Let $X \in L^{s}(\Omega,\mathcal{G},P)$. Note that,
$$
\int_{\Omega} |X|^{s}\,dP = \int_{\Omega} |X^r|^{s/r}\,dP.
$$ 

Hence, $|X|^{r} \in L^{s/r}(\Omega,\mathcal{G},P)$ and
$$
\Vert |X|^{r} \Vert_{s/r} = \left(\int_{\Omega} |X|^{r(s/r)}\,dP \right)^{r/s} = \left(\int_{\Omega} |X|^{s}\,dP \right)^{r/s} =\Vert X \Vert_{s}^{r}.
$$

By other hand, the function $1_{\Omega} \in L^{s/(s-r)}(\Omega,\mathcal{G},P)$. By H\"older-Riesz inequality (Theorem \ref{HRp}) we have that $1_{\Omega}|X|^{r}  = |X|^{r} \in L^{1}(\Omega,\mathcal{G},P)$ and thus
$$
\int_{\Omega} |X|^{r} \,dP \leq \Vert 1_{\Omega} \Vert_{s/(s-r)} \Vert |X|^{r} \Vert_{s/r}  = \Vert X \Vert_{s}^{r}.
$$

Consequently, $X \in L^{r}(\Omega,\mathcal{G},P)$ and $\Vert X \Vert_{r} \leq  \Vert X \Vert_{s}$.
\end{proof}

The Riesz-Fischer theorem \cite[Theorem 4.8]{Brezis(2011)} guarantees that the space $L^{p}(\Omega,\mathcal{G},P)$ with the norm $\Vert \cdot \Vert_{p}$ is a Banach space for all $p \in (1,\infty)$.

If we note, the identity that allowed us to conclude that Hilbert spaces are reflexive through the Milman-Pettis theorem (see Theorem \ref{MilmanPettis}) was the parallelogram identity, and in fact, it was sufficient to have an inequality. In $L^{p}$ spaces, there are certain inequalities that similarly allow us to conclude that they are reflexive spaces. These inequalities are known as the first and second Clarkson inequalities, and for their proof, we suggest consulting \cite{Rama(1978)} and \cite[Theorem 4.10]{Brezis(2011)}.

\begin{theorem}[Clarkson inequalities] \label{ClarksonInequality}
Let $p \in (1,\infty)$ and $X,Y \in L^{p}(\Omega,\mathcal{G},P)$.
\begin{itemize} 
                  \item[\textit{(1)}] If $p \in (1,2]$ then
$$
\left\Vert \frac{X+Y}{2}\right\Vert_{p}^{p/(p-1)} + \left\Vert \frac{X-Y}{2}\right\Vert_{p}^{p/(p-1)} \leq  \left( \frac{1}{2}\Vert X \Vert_{p}^{p} + \frac{1}{2}\Vert Y \Vert_{p}^{p}\right)^{1/(p-1)}.
$$
                  
                  \item[\textit{(2)}] If $p \in [2,\infty)$ then
                  $$
\left\Vert \frac{X+Y}{2}\right\Vert_{p}^{p} + \left\Vert \frac{X-Y}{2}\right\Vert_{p}^{p} \leq \frac{1}{2} \left( \Vert X \Vert_{p}^{p} + \Vert Y \Vert_{p}^{p}\right).
$$
              \end{itemize} 
\end{theorem}

Note that if $p=2$, the Clarkson inequalities imply the parallelogram identity in $L^{2}(\Omega,\mathcal{G},P)$.

\begin{corollary}
    $L^{p}(\Omega,\mathcal{G},P)$ is reflexive for all $p \in (1,\infty)$.
\end{corollary}

\begin{proof}
     Let $p \in (1,\infty)$. By Theorem \ref{MilmanPettis}, it is sufficient to prove that $L^{p}(\Omega,\mathcal{G},P)$ is uniformly convex (see Definition \ref{UniformConvexity}). Let $\varepsilon >0$ and let $X,Y \in L^{p}(\Omega,\mathcal{G},P)$ with $\Vert X\Vert_{p},\,\Vert Y \Vert_{p} \leq 1$ and $\Vert X-Y\Vert_{p} > \varepsilon$. 

{\scshape Case 1.} $p \in (1,2]$.

By Clarkson's first inequality we obtain
$$
\left\Vert \frac{X+Y}{2}\right\Vert_{p}^{q}  \leq  \left(\frac{1}{2}\, \Vert X \Vert_{p}^{p} + \frac{1}{2}\,\Vert Y\Vert_{p}^{p}\right)^{1/(p-1)} - \left\Vert \frac{X-Y}{2}\right\Vert_{p}^{q} < 1-\left(\frac{\varepsilon}{2}\right)^{q}
$$

Thus,
$$
\left\Vert \frac{X+Y}{2}\right\Vert_{p} < \left(1-\frac{\varepsilon^{p}}{2^{p}} \right)^{q} = 1-\left(1-  \left(1-\frac{\varepsilon^{q}}{2^{q}} \right)^{1/q} \right).
$$
and with $\delta:=\left(1-  \left(1-\tfrac{\varepsilon^{q}}{2^{q}} \right)^{1/q} \right)>0$ we conclude the result.

{\scshape Case 2.} $p \in [2,\infty)$. 

Applying Clarkson's second inequality to $X$ and $Y$, we obtain
$$
\left\Vert \frac{X+Y}{2}\right\Vert_{p}^{p}  \leq \frac{1}{2} \left( \Vert X \Vert_{p}^{p} + \Vert Y \Vert_{p}^{p}\right) -  \left\Vert \frac{X-Y}{2}\right\Vert_{p}^{p} < 1-\frac{\varepsilon^{p}}{2^{p}}
$$
from which it follows that
$$
\left\Vert \frac{X+Y}{2}\right\Vert_{p} < \left(1-\frac{\varepsilon^{p}}{2^{p}} \right)^{1/p} = 1-\left(1-  \left(1-\frac{\varepsilon^{p}}{2^{p}} \right)^{1/p} \right).
$$

With $\delta:=\left(1-  \left(1-\tfrac{\varepsilon^{p}}{2^{p}} \right)^{1/p} \right)>0$ the result holds.
\end{proof}

Therefore, using the Hahn-Banach theorem and the reflexivity of
$L^{p}$ spaces, we can obtain the Riesz representation theorem \cite[Theorem 4.11]{Brezis(2011)}.

\begin{theorem}[Riesz representation theorem on $L^p$] \label{Rieszp}
Let $p \in (1,\infty)$ and let $\Phi:L^{p}(\Omega,\mathcal{G},P) \to \mathbb{R}$ be a linear and continuous function. Then there exists a unique function $\psi \in L^{q}(\Omega,\mathcal{G},P)$ such that
$$
\Phi(X)=\int_{\Omega} \psi\cdot X\,dP \,\,\,\,\,\,\forall\,X \in L^{p}(\Omega,\mathcal{G},P).
$$
\end{theorem}

If $p=2$, then $q=2$ and in this case, Theorem \ref{Rieszp} turns out to be Theorem \ref{FrechetRiesz} in the Hilbert space $L^{2}(\Omega,\mathcal{G},P)$.

As in Section \ref{L2}, given $A \in \mathcal{F}$, we consider the function $T:L^{p}(\Omega,\mathcal{G},P) \to \mathbb{R}$ defined by
\begin{equation}\label{FunctionalTp}
T(X):=\int_{\Omega} X\cdot 1_{A}\,dP.    
\end{equation}

This function is clearly linear. Using the H\"older-Riesz inequality (Theorem \ref{HRp}), we obtain $|T(X)|\leq \Vert X \Vert_{p}$ for all $X \in L^{p}(\Omega,\mathcal{F},P)$. That is, $T$ is a linear and continuous function (see Definition \ref{CFunction}). We know from the Riesz representation theorem (Theorem \ref{Rieszp}) that there exists a unique $\xi \in L^{q}(\Omega,\mathcal{G},P)$ such that
$$
T(X)=\int_{\Omega} \xi\cdot X\,dP \quad \forall X \in L^{p}(\Omega,\mathcal{G},P).
$$

In particular, $T(1_{B})=\int_{\Omega} \xi \cdot 1_{B}\,dP$ for all $B \in \mathcal{G}$ causes $1_{B} \in L^{p}(\Omega,\mathcal{G},P)$ for all $B \in \mathcal{G}$. Now, since $1<q$, we have $\xi \in L^{1}(\Omega,\mathcal{G},P)$ by Corollary \ref{CHRp}. Thus, $\xi = E(1_{A}\,|\,\mathcal{G})$.

\section{The variational problem} \label{VariationalProblem}

The purpose of this section is to establish the existence of the conditional expectation of an event $A$ given a $\sigma$-algebra $\mathcal{G}$ as the critical point of a functional defined in certain Lebesgue spaces. Before doing so, we state and prove some preliminary results.

\subsection{Differentiable functionals}

Let $V=(V,\Vert \cdot \Vert_V)$ be a Banach space over $\mathbb{R}$ and $\mathcal{O}$ be an open set of $V$.

\begin{definition} \label{FrechetDif}
    Let $F:\mathcal{O} \to \mathbb{R}$ be a function. $F$ is differentiable at the point $u_{0} \in \mathcal{O}$ if there exists $L \in V^{\ast}$ such that:
$$
\lim_{v \to 0_V} \frac{| F(u_0+v)-F(v)-L(v) |}{\Vert v \Vert_V}=0.$$

$L$ is called the derivative of $F$ at $u_0$ and is denoted by $L:=F'(u_0)$.

$F$ is differentiable in $\mathcal{O}$ if it is differentiable at every point $u \in \mathcal{O}$. The function $F':\mathcal{O} \to V^{\ast}$ given by
$$
u \mapsto F'(u)
$$
is called the (Fréchet) derivative of $F$ in $\mathcal{O}$. If the function $F'$ is continuous, we say that $F$ is of class $\mathcal{C}^{1}$ on $\mathcal{O}$.
\end{definition}

\begin{definition} \label{DirDif}
Let $F: \mathcal{O} \to \mathbb{R}$ be a function. If, for $u \in \mathcal{O}$ and $v \in V$, the limit
$$
\lim_{t \to 0} \frac{F(u+tv)-F(u)}{t}
$$
exists, then its value is called the derivative of $F$ at the point $u$ and in the direction $v$.

A function $F:\mathcal{O} \to \mathbb{R}$ is Gateaux differentiable at point $u_{0} \in \mathcal{O}$ if, for every $v \in V$, the derivative of $F$ at point $u$ and in direction $v$ exists and the function $\mathcal{G}F(u):V \to \mathbb{R}$ given by
$$
\mathcal{G}F(u)v:=\lim_{t \to 0} \frac{F(u+tv)-F(u)}{t}
$$
is an element on $V^{\ast}$.

$F:\mathcal{O} \to \mathbb{R}$ is Gateaux-differentiable on $\mathcal{O}$ if it is Gateaux differentiable at every point $u \in \mathcal{O}$. The function given by
$$
\mathcal{G}F:\mathcal{O} \to V^{\ast},\,\,\,\,\,\,\,\,u \mapsto \mathcal{G}F(u),
$$
is called the Gateaux derivative of $F$.
\end{definition}

The following result relates the previous concepts and its proof can be found in \cite[Proposition 3.2.15]{Drabek(2007)}

\begin{theorem} \label{FGDif}
$F:\mathcal{O} \to \mathbb{R}$ is of class $\mathcal{C}^1$ in $\mathcal{O}$ if and only if $F$ is Gateaux-differentiable in $\mathcal{O}$ and $\mathcal{G}F:\mathcal{O} \to V^{\ast}$ is continuous. In this case, $F'=\mathcal{G}F$.
\end{theorem}

\begin{definition} \label{SecondDiff}
Let $F:\mathcal{O} \to\mathbb{R}$ be a function.
\begin{itemize}
    \item[(a)] If $F$ has a Fréchet derivative and $F':\mathcal{O} \to V^{\ast}$ is of class $\mathcal{C}^{1}$ on $\mathcal{O}$ then $F$ is class $\mathcal{C}^2$. The derivative of $F'$ is called the second derivative of $F$ and is denoted by $F''$.
    \item[(b)]  If $F$ is of class $\mathcal{C}^1$ in $\mathcal{O}$ then $F$ has a second Gateaux derivative in $u_0 \in \mathcal{O}$ if there exists $L \in \mathcal{B}(V,V^{\ast})$ such that, for every $v,w \in V$
    $$
     \lim_{t \to 0}\frac{ F'(u_{0} + t v)w - F'(u_{0})v - L(tv)w }{t}=0.
    $$

     We denote a second Gateaux derivative of $F$ in $u_0$ by $L:=\mathcal{G}^2F(u_{0})$. 
     
     $F:\mathcal{O} \to \mathbb{R}$ is twice Gateaux-differentiable on $\mathcal{O}$ if it has a second Gateaux derivative at every point $u \in \mathcal{O}$. The function given by
$$
\mathcal{G}^2F:\mathcal{O} \to \mathcal{B}(V,V^{\ast}),\,\,\,\,\,\,\,\,u \mapsto \mathcal{G}^2F(u),
$$
is called the second Gateaux derivative of $F$.
\end{itemize}
\end{definition}

As before, it is possible to relate the second Fréchet derivative to the second Gateaux derivative of a function $F:\mathcal{O} \to \mathbb{R}$ \cite[Proposition 1.6]{Willem(1996)}.

\begin{theorem} \label{SFG}
    $F$ is of class $\mathcal{C}^2$ in $\mathcal{O}$ if and only if $F$ has a continuous second Gateaux derivative in $V$. In this case, $F''=\mathcal{G}^2F$.
\end{theorem}

\begin{remark} \label{GBil}
Let $F:\mathcal{O} \to \mathbb{R}$ be a function of class $\mathcal{C}^2$ in $\mathcal{O}$.

In this text, we will consider the second Gateaux derivative as a continuous bilinear form, that is, for any $u \in \mathcal{O}$,
$$
\mathcal{G}^2F(u):V \times V \to \mathbb{R}
$$
that is linear in each entry and is continuous in $V \times V$ \cite[Remark 3.2.29]{Drabek(2007)}. Hence, the second Gateaux derivative on $u \in \mathcal{O}$ is given by
$$
\mathcal{G}^{2}F(u_{0})(v,w):=\lim_{t \to 0} \frac{F'(u_{0}+tv)w - F'(u_{0})w}{t}.
$$
\end{remark}

The reality is that the concept of Gateaux derivative, along with the Theorems \ref{FGDif} and \ref{SFG}, provides us with a very useful criterion to verify the differentiability of a functional and to compute its derivative.

\begin{proposition}\label{DifLinear}
If $L \in V^{\ast}$ then $L$ is of class $\mathcal{C}^{2}$ in $V$, $L'=L$ and $L''=0$.
\end{proposition}

\begin{proof}
Let $u \in V$. We have $L(u+v)-L(u)-L(v)=0$ for all $v \in V$. Thus, $L'(u)=L$ for all $u \in V$. Given that $L'(u+v)=L$ for all $u,v \in V$ then $L'$ is a constant function and therefore continuous. Consequently, $L''=0$ and $L$ are of class $\mathcal{C}^2$.
\end{proof}

\begin{proposition} \label{DifNorm}
    Let $H=(H,\langle \cdot,\cdot \rangle,\Vert \cdot \Vert_{H})$ be a Hilbert space over $\mathbb{R}$. The function $F:H \to \mathbb{R}$ given by
    $$
     F(u)=\frac{1}2\Vert u \Vert_{H}^{2}
    $$
    is of class $\mathcal{C}^2$, $F'(u)v=\langle u,v \rangle$ and $F''(u)(v,w)=\langle v,w \rangle$ for all $v,w \in H$.
\end{proposition}

\begin{proof}
Let $u \in H$. For every $v \in H$ and $t \in \mathbb{R}$ we have that
$$
\begin{aligned}
F(u+tv)-F(u) &=\frac{1}{2}\left( \Vert u \Vert_{H}^2 - 2t\langle u,v\rangle + t^2\Vert v \Vert_{H}^2  - \Vert u \Vert_{H}^2 \right)= \\& t\langle u,v\rangle + \frac{t^2}{2}\Vert v \Vert_{H}^2.
\end{aligned}
$$

Hence, 
$$
\lim_{t \to 0} \frac{F(u+tv)-F(u)}{t} = \lim_{t \to 0} \left( \langle u,v \rangle - \frac{t}{2}\Vert v \Vert_H^2 \right) = \langle u,v\rangle \,\,\,\,\,\,\,\forall\,v \in H.
$$

That is, $\mathcal{G}F(u)v=\langle u,v \rangle$ for all $v \in H$ and it's clear that $\mathcal{G}F(u) \in H^{\ast}$ by Riesz–Fréchet representation theorem (Theorem \ref{FrechetRiesz}).

Let $(u_j)$ be a sequence on $H$ such that $\Vert u_j - u \Vert_H \to 0$. By Cauchy–Schwarz inequality,
$$
|\mathcal{G}F(u_j)v-\mathcal{G}F(u)v|=|\langle u_j - u, v \rangle| \leq \Vert u_j - u \Vert_H \Vert v \Vert_H \,\,\,\,\,\forall\, v \in H.
$$

Thus, $\mathcal{G}F$ is continuous and hence $F$ is of class $\mathcal{C}^1$ on $H$ (see Theorem \ref{FGDif}).

Now, for all $v,w \in H$ and $t \in \mathbb{R}$,
$$
\mathcal{G}F(u+tv)w-\mathcal{G}F(u)w=\langle u+tv,w\rangle - \langle u,w \rangle = t\langle v,w \rangle.
$$

Thus,
$$
\lim_{t \to 0} \frac{\mathcal{G}F(u+tv)w-\mathcal{G}F(u)w}{t} = \langle v,w \rangle \,\,\,\,\,\,\forall v,w \in H
$$
and consequently, $\mathcal{G}F^2(u)(v,w)=\langle v,w \rangle$ for all $v,w \in H$ and it is clear that it is bilinear continuous. Theorem \ref{SFG} implies that $F$ is of class $\mathcal{C}^2$ and $F''=\mathcal{G}^2F$.
\end{proof}

\begin{definition}\label{CriticalPoint}
    Let $F:\mathcal{O} \to \mathbb{R}$ be a function of class $\mathcal{C}^{1}$. $u_0$ is a critical point of $F$ in $\mathcal{O}$ if $F'(u_0)v=0$ for all $v \in V$.
\end{definition}

For example, if $u_0 \in \mathcal{O}$ is a minimum point of $F:\mathcal{O} \to \mathbb{R}$ then $u_0$ is a critical point \cite[Theorem 2.A]{Zeidler(1995)}

According to the established objective, the following result provides a sufficient condition for finding minimum points of certain functionals \cite[Theorem 2]{Gelfand(1963)}.

\begin{theorem} \label{SufMin}
Let $F:\mathcal{O} \to \mathbb{R}$ be a function of class $\mathcal{C}^{2}$. $F$ has a minimum point $u \in \mathcal{O}$ if $F'(u)v=0$ for all $v \in V$ and $F''(u)(v,v) > 0$ for all $v \in V \smallsetminus \{0_V\}$.
\end{theorem}

\subsection{Dirichlet's principle}

Let $(\Omega,\mathcal{F},P)$ be a probability space, $\mathcal{B}=\{B_1,\ldots,B_N\}$ a finite collection of elements from $\mathcal{F}$ that form a partition of $\Omega$ and $\mathcal{G}:=\sigma(\mathcal{B})$.

Fix $A \in \mathcal{F}$, and we consider $(\Omega,\mathcal{G},P)$ the probability space given in Section \ref{usefulfacts}. The following result will allow us to interpret the conditional expectation of $A$ given the $\sigma$-algebra $\mathcal{G}$ as the unique minimum of the functional energy. This result is essentially Dirichlet's principle.

The functional energy on $L^{2}(\Omega,\mathcal{G},P)$ is the function $J:L^{2}(\Omega,\mathcal{F},P) \to \mathbb{R}$ given by
\begin{equation} \label{FunEn}
    J(X):=\frac{1}{2}\int_{\Omega} |X|^2\,dP - \int_{\Omega} X\cdot 1_{A}\,dP.
\end{equation}

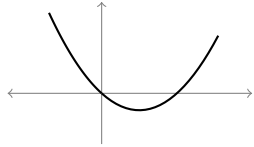
\begin{figure}[h!]
\centering
\begin{tikzpicture}[xscale=1,yscale=0.9]
	\draw[->,gray] (0,0) -- (2,0); \draw [->,gray] (0,0) -- (0,1.35); 
	\draw[<-,gray] (-1.25,0)--(0,0); \draw[-,gray] (0,0)--(0,-0.75);
\draw[domain= 0:1.55,thick] plot(\x,{\x^2 - \x} );
\draw[domain= -0.7:0,thick] plot(\x, {(\x)^2 - (\x)} );

\end{tikzpicture}
\caption{$J_{u}(t):=\frac{1}{2}\Vert tu \Vert_{2}^{2} - T(tu)$ for fixed $u \in L^{2}(\Omega,\mathcal{G},P)$ with $u \neq 0$ and $t \in \mathbb{R}$}
\end{figure}

If we prove that $J$ is of class $\mathcal{C}^1$ and we have that $\xi_0$ is a critical point of $J$ on $L^{2}(\Omega,\mathcal{G},P)$ then $\xi_0$ satisfies that
$$
J'(\xi_0)X=\langle \xi_0,X \rangle_{2} - T(X) =0\,\,\,\,\,\,\,\,\forall\,X \in L^{2}(\Omega,\mathcal{G},P).
$$

That is, we infer that the conditional expectation of $A$ given $\mathcal{G}$ is a critical point of the functional energy. In fact, let $\xi:=E(1_A\,|\,\mathcal{G}) \in L^2(\Omega,\mathcal{G},P)$. Then, $\int_{\Omega} \xi\cdot 1_{B} = T(1_B)$ for any $B \in \mathcal{G}$ and thus $\int_{\Omega} \xi \cdot Z = T(Z)$ for any simple random variable $Z \in L^2(\Omega,\mathcal{G},P)$ and, by consequence, $\int_{\Omega} \xi \cdot X = T(X)$ for any $X \in L^2(\Omega,\mathcal{G},P)$ because there exists a sequence of simple random variables $(Z_k)$ such that $Z_k \to X$ in $L^2(\Omega,\mathcal{G},P)$ (see \cite[Thereom 17.2 and Theorem 17.4]{Jacod(2004)}). We have the following result.

\begin{theorem}[Dirichlet's principle]
$E(1_{A}\,|\,\mathcal{G})$ is the unique minimum of the function $J:L^{2}(\Omega,\mathcal{F},P) \to \mathbb{R}$ given in (\ref{FunEn}).
\end{theorem}

\begin{proof}
Note that $J(X)=\frac{1}{2}\Vert X \Vert_{2}^{2} - T(X)$ with $T:L^{2}(\Omega,\mathcal{G},P) \to \mathbb{R}$ the linear and continuous function given in (\ref{FunctionalT2}). Hence, $J$ is of class $\mathcal{C}^2$ and and for $X \in L^{2}(\Omega,\mathcal{F},P)$ (see Propositions \ref{DifLinear} and \ref{DifNorm})
$$
J'(X)Y=\langle X,Y \rangle_{2} - T(Y) \,\,\,\,\,\,\,\,\forall\,Y \in L^{2}(\Omega,\mathcal{G},P).
$$

By Riesz–Fréchet representation theorem (Theorem \ref{FrechetRiesz}), $\xi:=E(1_{A}\,|\,\mathcal{G})$ exists and is the unique element in $L^{2}(\Omega,\mathcal{G},P)$ such that
$$
\langle \xi,Y\rangle_{2}-T(Y)=0 \,\,\,\,\,\,\,\,\forall\,Y \in L^{2}(\Omega,\mathcal{G},P).
$$

Thus, $J'(\xi)Y=0$ for all $Y \in L^{2}(\Omega,\mathcal{G},P)$.

Now, $J''(\xi)(Y,Y) = \langle Y,Y \rangle_{2} =E(Y^2) >0$ for all $Y \in L^{2}(\Omega,\mathcal{G},P)$ with $Y \neq 0$ (see Proposition \ref{DifNorm}). Theorem \ref{SufMin} states that $\xi=E(1_{A}\,|\,\mathcal{G})$ is a minimum of $J$ in $L^{2}(\Omega,\mathcal{G},P)$.
\end{proof}

The objective now is to extend the Dirichlet principle to $L^{p}$ spaces for $p \in (1,\infty)$, that is, we seek a functional $F:L^{p}(\Omega,\mathcal{G},P) \to \mathbb{R}$ of class $\mathcal{C}^2$ such that $F'(E(1_{A}\,|\,\mathcal{G}))Y=T(Y)$ for all $Y \in L^{p}(\Omega,\mathcal{G},P)$ with $T$ being the functional given in (\ref{FunctionalTp}).

Note that if $p \in (1,2)$ then exponential conjugate $q \in (2,\infty)$ and by Corollary \ref{CHRp} the functional $J_{q}:L^{q}(\Omega,\mathcal{G},P) \to \mathbb{R}$ given by
\begin{equation} \label{FunEnq}
J_{q}(X):=\frac{1}{2}\int_{\Omega} |X|^2\,dP - T(X)    
\end{equation}
is well-defined since $L^{q}(\Omega,\mathcal{G},P) \subset L^{2}(\Omega,\mathcal{G},P) \subset L^{p}(\Omega,\mathcal{G},P) \subset L^{1}(\Omega,\mathcal{G},P)$ because $1 < p < 2 < q < \infty$. We must note that we are considering the restriction of $T$ to $L^{q} (\Omega,\mathcal{G},P)$ which is clearly linear and continuous.

\begin{theorem} \label{DifLq}
    If $p \in (1,2)$, then the functional $J_{q}$ defined in (\ref{FunEnq}) is of class $\mathcal{C}^2$ on $L^{q}(\Omega,\mathcal{G},P)$, $J_{q}'(X)Y=\int_{\Omega} X\cdot Y\,dP - T(Y)$ and $J_{q}''(X)(Y,Z)=\int_{\Omega} Y\cdot Z\,dP$ for all $X,Y,Z\in L^{q}(\Omega,\mathcal{F},P)$.
\end{theorem}

\begin{proof}
Let $X \in L^{q}(\Omega,\mathcal{G},P)$. For every $Y \in L^{q}(\Omega,\mathcal{G},P)$ and $t \in \mathbb{R}$,
$$
\begin{aligned}
 J_q(X + tY) - & J_q(X) = \frac{1}{2}\int_{\Omega} (X+tY)^2\,dP -  T(X +tY) - 
 \\&\frac{1}{2}\int_{\Omega} |X|^2\,dP + T(X)=t\int_{\Omega} X\cdot Y\,dP - t T(Y) = \\& t\left(\int_{\Omega} X\cdot Y\,dP -T(Y) \right).
\end{aligned}
$$

Thus,
$$
\lim_{t \to 0} \frac{J_q(X + tY) - J_q(X)}{t} = \int_{\Omega} X\cdot Y\,dP - T(Y) \quad \forall \,Y \in L^{q}(\Omega,\mathcal{G},P)
$$
and $\mathcal{G}J_{q}(X)Y=\int_{\Omega} X\cdot Y\,dP - T(Y)$ for all $Y \in L^{q}(\Omega,\mathcal{G},P)$.

If $(X_j)$ is a sequence on $L^{q}(\Omega,\mathcal{G},P)$ such that $X_j \to X$ on $L^{q}(\Omega,\mathcal{G},P)$ by H\"older-Riesz inequality
$$
\left|\mathcal{G}J_{q}(X_j)Y-\mathcal{G}J_{q}(X)Y \right| \leq \int_{\Omega} |X_j-X||Y|\,dP \leq \Vert X_j - X\Vert_{q}\Vert Y\Vert_{p}.
$$

Hence, $\mathcal{G}J_{q}$ is continuous. Theorem \ref{FGDif} implies that $J_{q}$ is of class $\mathcal{C}^1$ and $J_{q}'=\mathcal{G}J_{q}$.

Now, for every $Y,Z \in L^{q}(\Omega,\mathcal{G},P)$ and $t \in \mathbb{R}$,
$$
\begin{aligned}
&\mathcal{G}J_{q}(X + tY)Z-\mathcal{G}J_{q}(X)Z =\int_{\Omega} (X+tY)\cdot Z\,dP - T(Z) - \\& \int_{\Omega} X\cdot Z\,dP + T(Z)=t\int_{\Omega} Y\cdot Z\,dP.
\end{aligned}
$$

Consequently,
$$
\begin{aligned}
\mathcal{G}^{2}J_{q}(X)(Y,Z)&= \lim_{t \to 0} \frac{\mathcal{G}J_{q}(X + tY)Z-\mathcal{G}J_{q}(X)Z}{t} \\&= \int_{\Omega} Y\cdot Z\,dP\quad\forall \,Y,Z \in L^{q}(\Omega,\mathcal{G},P).
\end{aligned}
$$

We can conclude that $J_{q}$ is of class $\mathcal{C}^2$ and $J_{q}''=\mathcal{G}J_{q}^2$ by Theorem \ref{SFG}.
\end{proof}

\begin{corollary}
    If $p \in (1,2)$ then $E(1_{A}\,|\,\mathcal{G})$ is the unique minimum of the function $J_{q}:L^{q}(\Omega,\mathcal{G},P) \to \mathbb{R}$ given by
$$
J_{q}(X):=\frac{1}{2}\int_{\Omega} |X|^2\,dP - \int_{\Omega} X \cdot 1_{A}\,dP
$$   
in $L^{q}(\Omega,\mathcal{G},P)$.
\end{corollary}

\begin{proof}
By Riesz representation theorem (Theorem \ref{Rieszp}), $\xi=E(1_{A}\,|\,\mathcal{G})$ is the unique element on $L^{q}(\Omega,\mathcal{G},P)$ such that
$$
T(X)=\int_{\Omega} \xi\cdot X \,dP \,\,\,\,\,\,\,\,\forall\,X \in L^{p}(\Omega,\mathcal{G},P).
$$

In particular, $T(Y)=\int_{\Omega} \xi \cdot Y\,dP$ for all $Y \in L^{q}(\Omega,\mathcal{G},P)$. 

Thus, by Theorem \ref{DifLq} we have $J_{q}'(\xi)Y=\int_{\Omega}\xi\cdot Y\,dP - T(Y) =0$ for all $Y \in L^{q}(\Omega,\mathcal{G},P)$ and $J_{q}''(\xi)(Y,Y)=\int_{\Omega} Y^2 \,dP >0$ for all $Y \in L^{q}(\Omega,\mathcal{G},P)$ with $Y \neq 0$. The result follows from Theorem \ref{SufMin}.
\end{proof}

\section{Law of Total Probability Proof} \label{TP}

Let $(\Omega,\mathcal{F},P)$ be a probability space, $\mathcal{B}=\{B_1,\ldots,B_N\}$ a finite collection of elements from $\mathcal{F}$ that form a partition of $\Omega$ and $\mathcal{G}:=\sigma(\mathcal{B})$.

In the previous sections, we have proven the existence and uniqueness of the conditional expectation of an event $A \in \mathcal{F}$ given $\sigma$-algebra $\mathcal{G}$. We conclude this work by providing a proof of the total probability theorem using the tools developed.

\begin{proof}[Proof of Theorem \ref{TotalProbability}]
Let $A \in \mathcal{F}$, and we consider $E(1_{A}\,|\,\mathcal{G})$. Given that $E(1_{A}\,|\,\mathcal{G})$ is $\mathcal{G}$ -measurable, by Proposition \ref{GMeasurable}, there exist constants $\alpha_{1},\ldots,\alpha_{N}$ such that
$$
E(1_{A}\,|\,\mathcal{G})=\sum_{j=1}^{N} \alpha_{j}1_{B_j}.
$$

Applying property \textit{(ii)} of the Definition \ref{ConditionalExpectation} to the previous expression, we obtain that
$$
P(A \cap B_{j})=\int_{B_{j}} 1_{A}\,dP = \int_{B_{j}} E(1_{A}\,|\,\mathcal{G})\,dP = \alpha_{j}P(B_{j}) \,\,\,\,\,\,\forall\,j = 1,\ldots,N.
$$

Thus, $E(1_{A}\,|\,\mathcal{G})=\sum_{j=1}^{N} P(A\,|\,B_{j})1_{B_j}$ and now applying property \textit{(iii)} of the Definition \ref{ConditionalExpectation} we conclude that
$$
P(A)=\int_{\Omega} 1_{A}\,dP = \int_{\Omega} E(1_{A}\,|\,\mathcal{G})\,dP = \sum_{j=1}^{N} P(A\,|\,B_{j})P(B_{j}).
$$

That is, we proved the total probability theorem using the Riesz representation theorems (Theorems \ref{FrechetRiesz} and \ref{Rieszp}).
\end{proof}

\bibliographystyle{amsplain}

\end{document}